\documentclass{birkjour}

\newcommand{\C}{{\mathbb C}}
\newcommand{\D}{{\mathbb D}}

\newcommand{\T}{{\mathbb T}}

\newcommand{\N}{{\mathbb N}}

\newcommand{\ov}{\overline}

\renewcommand{\th}{\theta}

\title[Two Problems on Coinvariant Subspaces]
{Two Problems on Coinvariant Subspaces\\
of the Shift Operator}
\author[K. M. Dyakonov]{Konstantin M. Dyakonov}

\address{%
ICREA and Universitat de Barcelona\\ 
Departament de Matem\`atica Aplicada i An\`alisi\\ 
Gran Via, 585\\ 
E-08007 Barcelona\\ 
Spain}

\email{konstantin.dyakonov@icrea.cat}
\thanks{Supported in part by grant MTM2011-27932-C02-01 from El Ministerio de Ciencia 
e Innovaci\'on (Spain) and grant 2009-SGR-1303 from AGAUR (Generalitat de Catalunya).}

\begin{document}
\begin{abstract}
Two problems are posed that involve the star-invariant subspace $K^p_\theta$ (in the Hardy space $H^p$) 
associated with an inner function $\theta$. One of these asks for a characterization of the extreme 
points of the unit ball in $K^\infty_\theta$, while the other concerns the Fermat equation 
$f^n+g^n=h^n$ in $K^p_\theta$. 
\end{abstract}

\maketitle

Let $H^p$ stand for the classical Hardy space on the disk $\D:=\{z\in\C:|z|<1\}$. As usual, 
we identify $H^p$-functions with their boundary values and treat $H^p$ as a subspace of $L^p(\T,m)$; 
here $\T:=\partial\D$ is the unit circle and $m$ is the normalized arclength measure on $\T$. 
The {\it shift operator} $S:H^p\to H^p$ acts by the rule $(Sf)(z)=zf(z)$, and we know from Beurling's 
theorem that the nontrivial $S$-invariant subspaces in $H^p$, with $p\in(0,\infty)$, are precisely those 
of the form $\th H^p$, where $\th$ is an inner function. (By definition, an {\it inner function} is an 
$H^\infty$-function whose modulus equals $1$ a.e. on $\T$.) For Beurling's theorem and other basic 
facts about $H^p$ spaces, see \cite[Chap.\,II]{G}. 

\par Now, the {\it $S$-coinvariant} (or {\it star-invariant}, or {\it model}) {\it subspace} generated by 
an inner function $\th$ is 
$$K^p_\th:=H^p\cap\th\,\ov{H^p_0},$$ 
where $H^p_0=zH^p$. Here and below, we restrict ourselves to the range $p\in[1,\infty]$, and we 
endow $K^p_\th$ with the usual $H^p$-norm $\|\cdot\|_p$. Such a subspace is invariant under the 
{\it backward shift operator} $S^*:f\mapsto(f-f(0))/z$ on $H^p$ and, for $p$ finite, it provides the 
general form of an $S^*$-invariant subspace. This last statement is easy to deduce, at least for 
$p\in(1,\infty)$, from Beurling's theorem by duality, since $K^p_\th$ is the annihilator (in $H^p$) of 
the $S$-invariant subspace $\th H^q$ (in $H^q$, with $p^{-1}+q^{-1}=1$) with respect to the standard 
pairing $\langle f,g\rangle=\int_\T f\bar g\,dm$. In particular, we have $K^2_\th=H^2\ominus\th H^2$. 
We finally observe that, for $p\in[1,\infty]$ and $\th$ inner, the subspace $K^p_\th$ coincides with 
the kernel of $T_{\bar\th}$, the Toeplitz operator with symbol $\bar\th$, in $H^p$. 

\par This note offers two open problems on $K^p_\th$ spaces that were previously posed by the author at 
the workshop \lq\lq Operator Theory and Harmonic Analysis" in Oberwolfach, Germany (November 2010) and at 
the CRM conference \lq\lq Invariant Subspaces of the Shift Operator" in Montreal, Canada (August 2013). 
\par When stating the first of these, and later on, we write $\text{\rm ball}(X)$ to denote the closed 
unit ball of a Banach space $X$. 

\medskip\noindent\textbf{Problem 1.} Given an inner function $\th$, characterize the extreme points 
of $\text{\rm ball}(K^\infty_\th)$. 

\medskip One may begin by recalling that the extreme points of $\text{\rm ball}(H^\infty)$ are precisely 
the unit-norm functions $f\in H^\infty$ with 
\begin{equation}\label{eqn:koshka}
\int_\T\log(1-|f|)\,dm=-\infty; 
\end{equation}
see \cite{dLR} or \cite[Chap.\,IV]{G}. It follows that every non-inner function 
in $\text{\rm ball}(K^\infty_\th)$ will be non-extreme for $\text{\rm ball}(H^\infty)$. 
Indeed, for $f\in K^\infty_\th$, we have 
\begin{equation}\label{eqn:sobaka}
\left(1-|f|^2\right)\th=\th-f\cdot\bar f\th\in H^\infty,
\end{equation}
because $\bar f\th\in H^\infty$. Assuming that $|f|\le1$ and $|f|\not\equiv1$ on $\T$, we see 
from \eqref{eqn:sobaka} that $1-|f|^2$ agrees with the modulus of a non-null $H^\infty$-function. 
Therefore, $\log(1-|f|^2)\in L^1(\T)$ and the integral in \eqref{eqn:koshka} is convergent. 
\par On the other hand, the space $K^\infty_\th$ (and its unit ball) will contain inner functions if and 
only if $\th(0)=0$. Those inner functions are then precisely the divisors of $\th/z$, and they are sure 
to be extreme points for $\text{\rm ball}(K^\infty_\th)$. It is the non-inner extreme points of this ball 
that we are concerned with. 
\par In \cite{DMRL}, Problem 1 was solved in the simplest case where $\th(z)=z^{N+1}$, with $N\in\N$. 
In this case, $K^\infty_\th$ becomes the space of polynomials of degree at most $N$, endowed with the 
sup-norm; we shall denote the latter space by $\mathcal P_N$. The solution given in \cite{DMRL} is, 
however, less elementary than one might at first expect. Namely, a non-inner (or equivalently, non-monomial) 
unit-norm polynomial $P\in\mathcal P_N$ is shown to be extreme for $\text{\rm ball}(\mathcal P_N)$ if 
and only if a certain Wronski-type matrix built from $P$ has maximal rank. It turns out that this 
rank condition cannot be rephrased, in any reasonably explicit way, in terms of the \lq\lq smallness" 
or the zeros of $1-|P|^2$. In fact, a simple construction from \cite{DMRL} produces two unit-norm 
polynomials, $P$ and $Q$, in $\mathcal P_N$ satisfying $1-|P|^2=2(1-|Q|^2)$ and such that $P$ is 
a non-extreme point of $\text{\rm ball}(\mathcal P_N)$, while $Q$ is extreme. 
\par The case of $\mathcal P_N$, let alone $K^\infty_\th$ with a general $\th$, therefore 
exhibits a higher level of complexity than $H^\infty$ (where the extreme points are described 
by \eqref{eqn:koshka}) or the space of {\it real-valued} trigonometric polynomials $P$ of degree $\le N$ 
(where the description involves only the number of zeros of $1-|P|^2$ on $\T$, see \cite{DMRL, KR, R}). 
When $\th$ is a finite Blaschke product, the extreme points of $\text{\rm ball}(K^\infty_\th)$ are 
probably describable in the spirit of \cite{DMRL}, but the general case seems to call for a new technique. 
\par Before moving on to the next problem, we mention that the extreme points of $\text{\rm ball}(K^1_\th)$ 
were characterized in \cite{DPAMS}. A more general result was actually proved there, the space 
in question being the kernel (in $H^1$) of an arbitrary Toeplitz operator $T_\varphi$ with 
$\varphi\in L^\infty$. When specialized to the case $\varphi=\bar\th$, with $\th$ inner, 
the result tells us that a unit-norm function $f\in K^1_\th$ is an extreme point of 
$\text{\rm ball}(K^1_\th)$ if and only if the inner factors of $f$ and $\bar z\bar f\th$ are 
relatively prime. This can be viewed as an extension of de Leeuw and Rudin's theorem that identifies 
the extreme points of $\text{\rm ball}(H^1)$ as outer functions; see \cite{dLR} or \cite[Chap.\,IV]{G}. 

\medskip\noindent\textbf{Problem 2.} Prove or disprove \lq\lq Fermat's last theorem" for $K^p_\th$: whenever 
$\th$ is inner, $p\in[1,\infty]$ and $n\in\{3,4,5,\dots\}$, no nontrivial triple of functions $f,g,h\in K^p_\th$ 
satisfies $f^n+g^n=h^n$. 

\medskip Here, by calling a triple $(f,g,h)$ nontrivial we mean that the three functions do not lie in a single 
one-dimensional subspace of $K^p_\th$. 
\par In the polynomial case (i.e., for $\th(z)=z^{N+1}$), Fermat's last theorem is known to be true: the 
equation $P^n+Q^n=R^n$ has no nontrivial polynomial solutions when $n\ge3$. (Such solutions do exist for 
$n=2$, as simple examples show.) It follows, then, that a similar result holds with rational functions in 
place of polynomials, and this settles Problem 2 in the case where $\th$ is a finite Blaschke product. 
The polynomial Fermat theorem, as stated above, is best proved as a consequence of the Mason--Stothers 
\lq\lq $abc$ theorem" (the prototype of the famous \lq\lq $abc$ conjecture" in number theory); 
the deduction consists in a three-line argument that can be found in \cite{DMA, GH, SS}. 
Now, the Mason--Stothers theorem provides a lower bound on the number of distinct zeros of 
the polynomial $PQR$ in terms of the degrees of $P$, $Q$ and $R$. In \cite{DMA}, we obtained 
some estimates of a similar nature for general holomorphic functions on the disk, while the 
case of entire functions was treated earlier in \cite{GH}. It seems that a suitable $K^p_\th$ version 
of the Mason--Stothers theorem, if available, might give us a clue to Problem 2.


\begin{thebibliography}{1} 

\bibitem{dLR} K. de Leeuw and W. Rudin, \textit{Extreme points and extremum problems in $H_1$.} 
Pacific J. Math. \textbf{8} (1958), 467--485. 

\bibitem{DPAMS} K. M. Dyakonov, \textit{Interpolating functions of minimal norm, star-invariant 
subspaces, and kernels of Toeplitz operators.} Proc. Amer. Math. Soc. \textbf{116} (1992), 1007--1013. 

\bibitem{DMRL} K. M. Dyakonov, \textit{Extreme points in spaces of polynomials.} Math. Res. Lett. 
\textbf{10} (2003), 717--728. 

\bibitem{DMA} K. M. Dyakonov, \textit{Zeros of analytic functions, with or without multiplicities.} 
Math. Ann. \textbf{352} (2012), 625--641. 

\bibitem{G} J. B. Garnett, \textit{Bounded analytic functions, Revised first edition.} Springer, 
New York, 2007. 

\bibitem{GH} G. G. Gundersen and W. K. Hayman, \textit{The strength of Cartan's version 
of Nevanlinna theory.} Bull. London Math. Soc. \textbf{36} (2004), 433--454. 

\bibitem{KR} A. G. Konheim and T. J. Rivlin, \textit{Extreme points of the unit ball in a space of 
real polynomials.} Amer. Math. Monthly \textbf{73} (1966), 505--507. 

\bibitem{R} H.-J. Rack, \textit{Extreme Punkte in der Einheitskugel des Vektorraumes der trigonometrischen 
Polynome.} Elem. Math. \textbf{37} (1982), 164--165. 

\bibitem{SS} T. Sheil-Small, \textit{Complex polynomials.} Cambridge Studies in Advanced Mathematics, 
75, Cambridge University Press, Cambridge, 2002. 

\end{thebibliography}
\end{document}